\documentclass[11pt]{article}

\usepackage[all]{xy}

\usepackage{amsmath}

\usepackage{amssymb}

\usepackage{bbm}

\newtheorem{theorem}{Theorem}[section]
\newtheorem{lemma}[theorem]{Lemma}
\newtheorem{corollary}[theorem]{Corollary}
\newtheorem{proposition}[theorem]{Proposition}
\newtheorem{exa}[theorem]{Example}

\newtheorem{exas}[theorem]{Examples}

\newtheorem{prope}[theorem]{Property}

\newtheorem{defini}[theorem]{Definition}

\newtheorem{rema}[theorem]{Remark}

\def\Z{\mathbbm{Z}}

\def\R{\mathbbm{R}}

\def\N{\mathbbm{N}}

\def\Q{\mathbbm{Q}}

\def\C{\mathbbm{C}}

\makeindex

\begin{document}

\date{Appeared in International Journal of Mathematics {\bf 20} no. 8
  (2009), 1057-1068.} 

\title{A
FINITENESS THEOREM FOR DUAL GRAPHS OF SURFACE SINGULARITIES} 
\maketitle

\begin{center}
PATRICK POPESCU-PAMPU\\
Universit\'e Paris 7 Denis Diderot \\
Institut de Math\'ematiques-UMR CNRS 7586\\
{\'e}quipe "G{\'e}om{\'e}trie et dynamique"\\
Site Chevaleret, Case  7012\\
75205 Paris Cedex 13\\
France.\\
e-mail: ppopescu@math.jussieu.fr
\end{center}

\begin{center}
JOS{\'E} SEADE\\
Unidad Cuernavaca del Instituto de Matem{\'a}ticas\\
Universidad Nacional Aut{\'o}noma de M{\'e}xico\\
Cuernavaca, M{\'e}xico\\
e-mail: jseade@matcuer.unam.mx
\end{center}

\begin{abstract}
Consider a fixed connected, finite graph $\Gamma$ and equip its
  vertices with
  weights $p_i$ which are non-negative integers. We show that
  there is  a finite number
of possibilities for the coefficients of the  canonical cycle of a
numerically Gorenstein surface singularity having $\Gamma$ as the
dual graph of  the minimal resolution, the weights $p_i$ of the
vertices being the arithmetic genera of the corresponding
irreducible components. As a consequence we get that if $\Gamma$
is not a cycle, then there is a finite number of possibilities of
self-intersection numbers which  one can attach to the vertices
which are of valency $\geq 2$, such that one gets the dual graph
of the minimal resolution of a numerically Gorenstein surface
singularity. Moreover, we describe precisely the situations when
there exists an infinite number of possibilities for the
self-intersections of the component corresponding to some fixed
vertex of $\Gamma$.
\footnote{ {Mathematics Subject Classification 2000: 14B05, 32S25,
  32S45}
\newline{  Key words and phrases: Gorenstein surface singularities,
  numerically Gorenstein, 
  canonical cycle,  dual graphs.}}
\end{abstract}


\section{Introduction}

Let $(X,0)$ be a germ of a normal complex analytic surface and
$(\tilde{X}, E) \stackrel{\pi}{\rightarrow} (X,0)$ a resolution of it,
its exceptional divisor $E$ \emph{not} being supposed to have
normal crossings. One associates to this resolution
a dual graph $\Gamma$, whose vertices are weighted by the absolute
values $e_i:=-E_i^2 >0$ of the self-intersections  and by the arithmetic
genera $p_i:=p_a(E_i)\geq 0$ of the corresponding irreducible
components $E_i$ of $E$. We see the two
collections of weights as functions with values in $\N$ defined on the
set of vertices of $\Gamma$, and we denote them by $\underline{e}$ and
$\underline{p}$ respectively.

 It is known since the works of Du Val \cite{V 44} and Mumford \cite{M
   61} that the
 intersection form associated to the weighted  graph $(\Gamma,
 \underline{e})$ is  negative definite. In this case we say that
 $(\Gamma, \underline{e})$ and $(\Gamma, \underline{e},
 \underline{p})$ are also \emph{negative definite}.
 Conversely,  Grauert  \cite{G 62} showed that if a compact connected
 reduced divisor on a smooth complex analytic surface has a weighted
 dual graph which satisfies this
 condition, then it may be contracted
 to a normal singular point of an analytic surface.

An important problem, studied by several authors  (see for
 instance Yau \cite{Y 79} and Laufer \cite{L 87}), is \emph{to decide,
 amongst negative definite dual
 graphs $(\Gamma, \underline{e}, \underline{p})$, which ones
 correspond to isolated hypersurface singularities in  $\C^3$}.
This  problem is rather deep and has resisted several attempts to solve it.
This article arose as a new attempt to give a step in that
direction. To be precise,
we look at the following weaker problem:

\medskip
 \emph{Describe the dual graphs corresponding to numerically
   Gorenstein normal surface singularities.}
\smallskip

\noindent This second problem was also posed by Jia, Luk \& Yau \cite{JLY 07}.

Recall that the normal germ $(X,0)$ is called \emph{Gorenstein} if the
canonical line bundle is
holomorphically trivial on a pointed neighborhood of $0$ in $X$. It is called
\emph{numerically Gorenstein} if the same line bundle is smoothly
trivial. An isolated
hypersurface singularity is a particular case of isolated complete
intersection singularity,
which is a particular case of a Gorenstein isolated singularity, which
is a particular case
of numerically Gorenstein isolated singularity.  This explains in
what sense the previous
problem is weaker than the initial one.

This article is devoted to  the study of the dual  graphs of
numerically Gorenstein surface singularities; for short, we call
these singularities \emph{n-Gorenstein}. The starting point is
Durfee's observation in  \cite {D 78} (see also    \cite{LS 95}),
that an isolated surface singularity  $(X,0)$ is n-Gorenstein
iff the \emph{canonical cycle} $Z_{can}$ of every  resolution
is integral. This 2-cycle, supported on the exceptional divisor,
is uniquely characterized by the fact that it satisfies the
adjunction formula:
$$ 2p_a(E_i)-2 = E_i^2 + Z_{can} \cdot E_i\,,$$ for each irreducible
component $E_i$
of the exceptional divisor, where $p_a(E_i)$ is the arithmetic
genus of the (possibly singular) curve $E_i$.
This system of equations shows that the coefficients  of
$Z_{can}$ are determined by the dual graph of the considered resolution,
decorated by the weights $e_i$ and $p_i$ attached to the
vertices. Therefore, we can speak about \emph{n-Gorenstein graphs}
$(\Gamma,  \underline{e}, \underline{p})$.

Given an arbitrary
finite, connected, unoriented graph $\Gamma$, whose set of vertices is
denoted $V(\Gamma)$, every choice of sufficiently large positive
weights $(e_i)_{i \in V(\Gamma)}$ makes it have negative definite
intersection form (see \cite{LS 95}). More precisely, there exist
weights $(e_i^0)_{i  \in V(\Gamma)}$
such that one has a negative definite intersection form whenever
$e_i\geq e_i^0 \: $  for all  $\: i \in V(\Gamma)$. By Grauert's theorem,
every choice of
arithmetic genus $p_i$ for each vertex turns $(\Gamma, \underline
e, \underline p)$  into the dual graph of a resolution of a normal surface
singularity. Furthermore (see \cite[Theorem 2.10]{LS 95}),
for each fixed function $\underline{e}$ such that $(\Gamma,
\underline{e})$ is negative defnite,  there
are infinitely many choices of a function $\underline{p}$ that make
$(\Gamma, \underline{e}, \underline{p})$  n-Gorenstein.

In this paper we  look at the converse situation,  which is  more
delicate:

\medskip
\emph{If we fix $(\Gamma, \underline{p})$,
  how many choices are there for $\underline{e}$  such
that $(\Gamma, \underline{e}, \underline{p})$ is a   n-Gorenstein graph?}
\medskip

As noticed before, the knowledge of $(\Gamma, \underline{e},
\underline{p})$ determines the coefficients of the canonical cycle of
a singularity whose minimal resolution has this dual graph.
We succeeded in proving the following finiteness result, which is a
consequence of the main theorem of this paper (see Theorem \ref{mainth}):

\medskip

\noindent \emph{Let $\Gamma$ be an arbitrary finite, connected, unoriented
graph equipped with  weights $p_i \ge 0$ assigned to each
vertex $i\in V(\Gamma)$.  Then there are at most a finite number of choices of
weights  corresponding to the coefficients of the
canonical cycle of the minimal resolution of a n-Gorenstein
surface singularity
with dual graph $\Gamma$ and arithmetic genera $p_i$.}
\medskip

As a consequence of this result we get that if $\Gamma$ is not a cycle,
then there is at most a finite number of choices of weights
$e_i$ for the vertices with valence at least $2$,
when one varies $\underline{e}$ such that $(\Gamma, \underline{e},
\underline{p})$ is n-Gorenstein and minimal (that is, without
vertices $i$ such that $p_i=0$ and $e_i =1$).  Moreover, we describe precisely
in all the cases the
non-finiteness appearing in the choice of the weights $(e_i)_{i
\in V(\Gamma)}$ (see Proposition \ref{measfin}).

Before describing briefly the content of each section, we would like
to mention that, as Jonathan Wahl told us,  it is unknown if each
n-Gorenstein graph occurs as the dual graph associated to a resolution
of a Gorenstein (and not merely n-Gorenstein) normal surface
singularity.

\medskip

In Sec. 2 we explain the necessary background about dual
graphs, associated quadratic forms, anticanonical cycles and
numerically Gorenstein singularities. In Sec. 3 we give
examples of families of numerically Gorenstein singularities,
whose study allowed us to conjecture the results proved in the
following sections. In Sec. 4 we prove our main theorem,
and in the last one we describe precisely the vertices $i$ of the graph
$\Gamma$, to which may be associated in an infinite number of ways a
value $e_i$, extendable to a weight $\underline{e}$, making $(\Gamma,
\underline{p}, \underline{e})$ n-Gorenstein. This places the
examples of the second section in a clearer light. We conclude with
some questions.

\section{Dual graphs and anticanonical cycles} \label{dual}

Let $(X,0)$ be a germ of normal complex analytic surface. Denote by
$(\tilde{X}, E)\stackrel{\pi}{\rightarrow} (X,0)$ its \emph{minimal
  resolution}, where
$E$ is the reduced fibre over $0$. Therefore $E$ can be regarded as a connected
reduced effective divisor in $\tilde{X}$, called \emph{the exceptional divisor}
of $\pi$. Note that $E$ has not necessarily normal
crossings. In particular, its irreducible components are not
necessarily smooth.

Denote by $\Gamma$ the \emph{dual (intersection) graph} of $E$: its vertices
correspond bijectively to the components of $E$ and between two distinct
vertices $i$ and $j$
there are as many (unoriented) edges as the intersection number
$e_{ij}:=E_i \cdot E_j\geq 0$ of the corresponding components. In particular,
$\Gamma$ has no loops. Moreover, each vertex
$i$ of $\Gamma$ is weighted by
the number $e_i$, where $-e_i:=E_i^2$ is the self-intersection number
of the associated component $E_i$ inside $\tilde{X}$.

Denote by $V(\Gamma)$ the set of vertices of $\Gamma$ and by
$\underline{e}\in \Z^{V(\Gamma)}$ the function which associates to each vertex
its weight.
To the weighted graph $(\Gamma, \underline{e})$ is
associated a canonical quadratic form on the real vector space
$\R^{V(\Gamma)}$,
called \emph{the intersection form} associated to the resolution
$\pi$:
\begin{equation} \label{formint}
 Q_{(\Gamma, \underline e)} (\underline{x}): = \sum_{i \in V(\Gamma)}
( -e_i x_i^2 +
   \sum_{\stackrel{j \in V(\Gamma)}{j\neq i}} e_{ij} x_i x_j)\;=\;
   \sum_{i \in V(\Gamma)} x_i ( -e_i x_i + \sum_{\stackrel{j \in
       V(\Gamma)}{j\neq i}} e_{ij}  x_j).
\end{equation}

  Often in the literature one fixes an order of the components of $E$, which
  makes $\underline{e}$
  appear as a finite sequence and $(e_{ij})_{i,j}$ appear as a matrix. We
  did not choose to do so in order to emphasize that there is no natural order
  and that our considerations do not depend on any such choice.

Du Val \cite{V 44} and Mumford \cite{M 61} proved that the
intersection form $Q_{(\Gamma, \underline e)}$
is negative definite. In particular, $e_i >0$  for all $i \in
V(\Gamma)$. Conversely, Grauert \cite{G 62} proved that if the
form associated to a reduced compact effective divisor $E$ on a
smooth surface is negative definite, then $E$  can be contracted
to a normal singular point of an analytic surface.
\medskip

For the following considerations on arithmetic genera, the adjunction
formula and
the anti-canonical cycle, we refer to Reid \cite{R 97} and Barth,
Hulek, Peters \& Van de Ven
\cite{BHPV 04}.

If $D$ is an effective divisor on $\tilde{X}$ supported on $E$, then it may be
interpreted as a (non-necessarily reduced) compact curve, with
associated structure sheaf $\mathcal{O}_D$. Its \emph{arithmetic genus}
$p_a(D)$ is by definition equal to $1-\chi(\mathcal{O}_D)$. It satisfies
\emph{the adjunction formula}:
\begin{equation} \label{adj}
   p_a(D):= 1 +\frac{1}{2}( D^2 + K_{\tilde{X}}\cdot D)
 \end{equation}
 where $K_{\tilde{X}}$ is any canonical divisor on $\tilde{X}$. This
 allows to extend the definition to \emph{any} divisor supported on $E$, not
 necessarily an effective one.

Denote by $p_i$ the arithmetic genus of the curve $E_i$  for all $
i \in V(\Gamma)$, and
by $g_i$ the arithmetic genus of its normalization, equal to its
topological genus.
Both genera are related by the following formula:
\begin{equation} \label{genus}
   p_i=g_i+ \sum_{P \in{E_i}} \delta_P(E_i)
 \end{equation}
where $\delta_P(E_i) \geq 0$ denotes the delta-invariant of the point $P$ of
$E_i$, equal to the
number of ordinary double points concentrated at $P$. One has
$\delta_P(E_i)>0$ if and
only if $P$ is singular on $E_i$. We deduce from (\ref{genus}) that:
\begin{equation} \label{carsp}
  p_i=0 \mbox{ if and only if }
  E_i \mbox{ is a smooth rational curve.}
\end{equation}

At this point, we have two weightings for the vertices of the graph
$\Gamma$: 
the collection $\underline{e}$ of self-intersections and the collection
$\underline{p}$ of arithmetic genera of the associated irreducible
components. If $E$ is a divisor with normal crossings and moreover all
its components are smooth, then the doubly weighted graph
$(\Gamma, \underline{e}, \underline{p})$ determines the embedded
topology of $E$ in $\tilde{X}$ (see Mumford \cite{M 61}). In general this
is not the case, because these numerical data do not determine the types
of singularities of $E$. Nevertheless, they determine them, and
consequently the embedded topology of $E$, up to a finite ambiguity. Indeed,
there are a finite number of embedded topological types of germs
of reduced plane curves $(C,c)$ having a given value $\delta_c(C)$
(see Wall \cite[page 151]{W 04}).

As the quadratic form $Q_{(\Gamma, \underline e)}$ is negative definite,
there exists
a unique divisor \emph{with rational coefficients} $Z_K$ supported on
$E$ such that:
\begin{equation} \label{numcycl}
    Z_K\cdot E_i=- K_{\tilde{X}}\cdot E_i\,, \: \;\hbox{for all} \,\: i
    \in V(\Gamma).
\end{equation}
We call $Z_K$ the \emph{anti-canonical cycle} of $E$ (or of the
resolution $\pi$). The name is motivated by the fact that whenever
$(X,0)$ is Gorenstein, $-Z_K$ is a canonical divisor on
$\tilde{X}$ in a neighborhood of $E$. With the notations of the introduction,
$Z_{can}=-Z_K$. The sign in the previous definition is motivated
by the following well-known result:

\begin{lemma} Assuming (as we do) that the resolution is minimal,
   $Z_K$ is an effective divisor.
\end{lemma}

\noindent{\bf Proof.} From (\ref{adj}) and (\ref{numcycl}) we get
    $Z_K\cdot E_i =
  -e_i-2p_i+2$. This number is necessarily non-positive. This is clear if
  $p_i\geq 1$. If instead $p_i=0$, by (\ref{carsp}) we see that $E_i$
  is a smooth
  rational curve. As $\pi$ is supposed to be the minimal resolution of
  $(X,0)$, we get $e_i\geq 2$ by Castelnuovo's criterion, which shows
  again that $-e_i-2p_i+2\leq0$.
  Therefore $Z_K\cdot E_i \leq 0$ for all $\: i\in V(\Gamma)$,
  which implies  that $Z_K$ is effective (cf. the proof of
 Proposition 2 in \cite{Artin}).  \hfill $\Box$
  \medskip

Denote: $$Z_K= \sum_{i\in V(\Gamma)}z_i E_i.$$ 
The previous lemma shows that
$\underline{z}\in \Q^{V(\Gamma)}_{\geq 0}$. By the adjunction formulae
(\ref{adj}) and
the relations (\ref{numcycl}), we
get the following system of equations relating $\underline e, \underline p$ and
$\underline z$:
\begin{equation} \label{eqnum}
   2p_i-2=-(z_i-1)e_i- \sum_{\stackrel{j \in V(\Gamma)}{j \neq i}} z_j e_{ij}
 \end{equation}

\begin{defini} The singularity $(X,0)$ is called \textbf{numerically
    Gorenstein} or
\textbf{n-Gorenstein} if $Z_K$ is an integral divisor.  As the
coefficients $\underline Z = \{z_i\}$ of $Z_K$ depend only on
 the decorated graph $(\Gamma,\underline p,\underline e)$, we also say
 that this graph is \textbf{n-Gorenstein}.
\end{defini}

Recall now that the \emph{Du Val singularities}, also known as
\emph{Kleinian singularities, rational double points or simple surface
singularities} (see Durfee \cite{D 79}) are, up to isomorphism,
the surface singularities of the form $\mathbb C^2/G$, where $G$
is a finite subgroup of $SU(2)$. For these singularities the
minimal resolution has $Z_K = 0$, the dual graph $\Gamma$ is one
of the trees $A_n, D_n, E_6, E_7, E_8$ and $p_i=0, e_i=2$ for all the vertices
$i$ of $\Gamma$.

\begin{lemma} If $(\Gamma,\underline p,\underline e)$ is not one of
  the Dynkin diagrams $\{A_n,D_n,E_6,E_7,E_8\}$ corresponding to the
  Du Val singularities
and is n-Gorenstein, then $z_i > 0$ for all $i\in V(\Gamma)$.
\end{lemma}

 \noindent{\bf Proof.} Suppose that $z_i = 0$. The
previous equation implies that $2p_i-2<0$, thus $p_i=0$. Therefore
(\ref{eqnum}) may be written:
$$-2 = -e_i - \sum_{\stackrel{j \in V(\Gamma)}{j\neq i}} z_j e_{ij}\,.$$
The hypothesis that the resolution is minimal shows that $e_i \ge
2$, since $p_i =0$. Hence $e_i=2$ and $z_j =0$ for all the
neighbors $j$ of $i$. Extending this argument step by step and using
the connectedness of $\Gamma$, one gets $z_j =
0$ and  $e_j = 2$ for all $j \in V(\Gamma)$. Therefore, the
decorated graph must be as stated, by a classical characterization
of Du Val singularities (see \cite[page 19]{BHPV 04}). \hfill
$\Box$
\medskip

In the sequel, we suppose that $(\Gamma,\underline p,\underline e)$ is not one
of the Dynkin diagrams $\{A_n,D_n,E_6,E_7,E_8\}$. By the previous
lemma, $z_i \geq 1$  for all $\: i \in V(\Gamma)$. Let us introduce new
variables, for simplicity:
\begin{equation} \label{change}
 \left\{ \begin{array}{l}
 n_i := z_i -1  \ge 0\,,\\
 v_i :  = \displaystyle{\sum_{\stackrel{j \in V(\Gamma)}{j\neq i}}
  }e_{ij}\geq 0, \\
 q_i : = v_i + 2 p_i -2 \ge -2.
 \end{array} \right.
\end{equation}

Then the adjunction formulae (\ref{eqnum}) become:

\begin{equation} \label{eqnumter}
  \{ e_i n_i \,=\, q_i + \sum_{\stackrel{j \in V(\Gamma)}{j\neq i}}
  \,e_{ij}n_j \}_{ i \in V(\Gamma) } .
 \end{equation}

 If $i\in V(\Gamma), \: v_i$ is the \emph{valency} of $i$, that is,
 the number of edges
 connecting it to other vertices.
 If $\Gamma$ is homeomorphic to a circle and $p_i=0$
 for all  $ i\in V(\Gamma)$,
 we say that $(\Gamma, \underline{p})$ is a \emph{cusp-graph}. The
 name comes from
 the fact that such graphs appear as dual resolution graphs of so-called
 \emph{cusp surface singularities} (see Brieskorn \cite[page 54]{B 00}
 or Looijenga \cite[page 16]{L 84}).

We summarize below the previous discussion:

\begin{proposition} Let $(X,0)$ be a n-Gorenstein surface
singularity which is not a Du Val singularity. Let $E =\sum E_i$
be the reduced exceptional divisor of its minimal resolution and
let $Z_K =\sum z_i E_i$ be the anti-canonical cycle,  characterized by
the adjunction formulae (\ref{eqnum}).
Then $Z_K-E$ is an effective divisor and, with the notations of 
(\ref{change}), the adjunction formulae become
the formulae (\ref{eqnumter}).
\end{proposition}

\section{Examples} \label{exampl}

In this section we present some of the examples which led us to
conjecture the results proved in the next two sections. In each one of
them, we fix a weighted graph $(\Gamma, \underline{p})$ and we look for
the weights $\underline{e}$ which make $(\Gamma, \underline{e},
\underline{p})$ correspond to a n-Gorenstein singularity.

i) Suppose that  the graph $\Gamma$ has only 1 vertex, no edges, and
that we equip the vertex
 with some weight $p \ge 0$.  Then the anti-canonical cycle is $Z_K =
 z E$ for some $z \in \N$,
 where $E$ represents a (possibly singular) irreducible projective curve of
 arithmetic genus $p$.  Denoting $e:=-E^2$, the adjunction formulae
 (\ref{eqnum})  imply:
$$ z =  \frac{2p-2}{e} + 1\,.$$
Hence we have the dual graph of an n-Gorenstein singularity iff
the weights $(p, e)$  are chosen so that
$\frac{2p-2}{e}$ is an integer.
Obviously, except when $p=1$, there are finitely many  choices of
such weights for a fixed $p$.

\medskip

ii) Consider  a graph $\Gamma$ with two vertices $1$ and $2$ and
one edge between them, and equip the vertices  with genera $p_1=1$ and
$p_2=2$ respectively.  The adjunction system (\ref{eqnumter}) becomes:
$$\left\{ \begin{array}{l}
   e_1 n_1 - n_2 =1 \\
    e_2 n_2 - n_1=3
    \end{array} \right.  .$$
An easy computation shows that there are exactly $8$ solutions
$(n_1, n_2; e_1, e_2)$ of the system, as follows: $(5,4; 1,2)$,
$(3,2;1,3)$, $(2,1;1,5)$, $(4,7;2,1)$, 
$(1,1; 2,4)$, $(2,5; 3,1)$, $(1,2; 3,2)$, $(1,4; 5,1)$.

\medskip

iii) Consider the quotient-conical singularities
of Dolgachev \cite{Do 74}: given any cocompact fuchsian group $G$
of signature $\{g;\alpha_1, ..., \alpha_n\}$, we may let it act on
$TH \cong H \times \C$ via the differential: $h \cdot (z,w)
\mapsto (h(z), h_*(z) \cdot w)$, where $H$ is the upper half plane
in $\C$. The surface $TH / G$ contains $H /  G$ as a divisor that
can be blown down analytically. The result is a normal surface
singularity $(X, 0)$ whose abstract boundary (or link) is
diffeomorphic to $PSL(2,\R) /G$ and which has a
resolution with dual graph a star with a center representing a curve
$E_0$ of genus $g$ and weight $e_0=2g- 2+ n$; it has  $n$ branches
of length 1, each with an end-vertex $i$
 that represents a curve of genus $0$ and weight $e_i=\alpha_i$. The
 $\alpha_i$ can take any
 values $\ge 2$. The anti-canonical cycle is $Z_K = 2 E_0 + \sum_{i=1}^n E_i$.
 Thus, given such a graph, equipped with the corresponding genera, there
 are infinitely many
 choices of weights $e_i$ for
 the vertices of valency 1 which make it correspond to
 n-Gorenstein
 singularities (cf.  \cite{D 83}, \cite{Ne 83}).

\medskip

iv) Consider as a final example the dual graph of a cusp
singularity (see \cite{B 00} or \cite{L 84}).
 This is a cycle of
finite length; all its vertices represent smooth rational curves
$E_1,...,E_n$ with $e_i \geq 2 $ for all $ i\in V(\Gamma)$,  and
at least one vertex $i$ satisfies  $e_i\geq 3$. The anti-canonical
cycle is $Z_K = E_1 + \cdots + E_n$, regardless of the weights
$\underline{e}$, which shows that such singularities are
n-Gorenstein. We see that all the choices of weights
$\underline{e}$ satisfying the previous inequality are good, for
every choice of genera. Therefore there exists an infinite number
of possibilities for each weight $e_i$.

\section{The main theorem} \label{main}

In this section we give a structure theorem about the set of solutions
of systems of the form (\ref{eqnumter}), where we drop the conditions
$q_i \geq -2$, which are necessarily satisfied (see (\ref{change})) if
the system corresponds to potential surface singularities.

If $i,j$ are distinct vertices of $\Gamma$, we denote by
$i\leftrightarrow j$
the fact that they are adjacent, that is, connected by at least
one edge.

\begin{theorem} \label{mainth}
Consider a graph $\Gamma$ decorated with weights   $\underline q \in
\Z ^{V(\Gamma)}$,
and the system of equations in the unknowns $ (\underline n ,
\underline e) \in
(\N)^{V(\Gamma)} \times (\N^*)^{V(\Gamma)} $:
  \begin{equation} \label{eqgen}
    \{e_i n_i = q_i + \sum_{\stackrel{j\in
        V(\Gamma)}{j\leftrightarrow i}} e_{ij}
    n_j \}_{i \in V(\Gamma)}.
  \end{equation}
  Then there exist at most finitely many weights $ \underline n$
  which can be extended
  to solutions $ (\underline n , \underline e)$ of the previous system,
such that the quadratic form $Q_{(\Gamma, \underline e)}$
is negative definite.
\end{theorem}

\noindent
{\bf Proof.} By working on examples, we noticed that we got
contradictions if we searched for  solutions of the system
(\ref{eqgen}) by traveling continuously on the graph $\Gamma$,
starting from a value $n_i$  which was too
big. As we were unable to find a precise description of what ``too
big'' meant, we had the idea to search a contradiction starting from a
sequence of solutions with unbounded values of $\underline{n}$. This
idea worked, as we explain now.

Suppose that there exists a
sequence of solutions  $ ((\underline n^{(k)} , \underline e^{(k)})
\in (\N)^{V(\Gamma)} \times (\N^*)^{V(\Gamma)} )_{k \geq 1}$ such
that:
$$N^{(k)}: = \hbox{max}_{i \in V(\Gamma)}\{n_i^{(k)}\}_{k \rightarrow_
  \infty} \rightarrow +\infty .$$

Selecting subsequences if necessary, we may assume that:

i)  there exists $i_o \in V(\Gamma)$ such that $n_{i_o}^{(k)} =
N^{(k)}$, for all $k\geq 1$;

ii) for all $i \in V(\Gamma)$, there exists $\displaystyle{\lim_{k \to
    \infty}} \frac{n_i^{(k)}}{n_{i_o}^{(k)}} =: \nu_i \in [0,1]$.

Set:  $$P := \{i \in V(\Gamma) \vert  \nu_i > 0 \};$$ then $P \ne
\emptyset$ (indeed, $i_o \in P$, as $\nu_{i_o} = 1$ by i)).

Let $\Gamma_P$ be the subgraph of $\Gamma$ spanned by $P$ (that is, the
subgraph of $\Gamma$ whose set of vertices is $P$ and whose edges are all
the edges of $\Gamma$ which connect two elements of $P$). It can be 
non-connected.
Denote by $\Gamma_{(P,i_o)}$  \emph{the connected component of $\Gamma_P$
which contains the vertex $i_o$}.

For all $i \in V(\Gamma_{(P,i_o)})$, one has by (\ref{eqgen}):

$$e_i^{(k)} \frac{n_i^{(k)}}{n_{i_o}^{(k)} } = \frac{q_i}{n_{i_o}^{(k)} }+
\sum_{\stackrel{j \in V(\Gamma)}{j\leftrightarrow i}} e_{ij}\,
\frac{n_j^{(k)}}{n_{i_o}^{(k)} } .$$

By assumption ii) and the
construction of $\Gamma_P$, the following limits exist:
$$\: \displaystyle{\lim_{k \to
    \infty}} \frac{n_i^{(k)}}{n_{i_o}^{(k)}} = \nu_i>0  \,,$$
  for all $i \in V(\Gamma_{(P,i_o)}) $,
    and:  $$  \: \displaystyle{\lim_{k \to \infty}}
\frac{n_j^{(k)}}{n_{i_o}^{(k)}} = 0\,,$$
for all the vertices $j \in V(\Gamma)\setminus V(\Gamma_{(P,i_o)})$ which are
connected to the vertex $i$.
Moreover, as $\displaystyle{\lim_{k \to \infty}}N^{(k)}=\infty$ and
 $N^{(k)} = n_{i_o}^{(k)}$  for all $\: k \geq 1 $ by assumption i),
 we see that
$ \displaystyle{\lim_{k \to \infty}} \frac{q_i}{n_{i_o}^{(k)} }=0$.

Thus one gets that for all $i \in V(\Gamma_{(P, i_o)})$,
$\displaystyle{\lim_{k \to \infty}}  e_i^{(k)}$ exists and: 

\begin{equation} \label{lime}
   \epsilon_i := \displaystyle{\lim_{k \to \infty}} e_i^{(k)}  =
   \frac{1}{\nu_i} \,
    \sum_{\stackrel{j \in V(\Gamma_{(P, i_o)})}{j\leftrightarrow i}
    }e_{ij} \nu_j < + \infty .
\end{equation}

As all the numbers $e_i^{(k)}$ are integers, this shows that  there
exists $k_o$ such that
for all $k \ge k_o$ and for all $i \in V(\Gamma_{P,i_o})$ one has:
$e_i^{(k)}=\epsilon_i$.
Therefore, by equation (\ref{lime}):

\begin{equation} \label{eqlim}
    e_i^{(k_o)} \, \nu_i = \sum_{\stackrel{j \in V(\Gamma_{(P,
      i_o)})}{j\leftrightarrow i} } e_{ij} \nu_j\,, \: \hbox {
   for all} \:\, i \in V(\Gamma_{(P, i_o)}).
\end{equation}

Define $\underline{\mu}\in \Z ^{V(\Gamma)}$ by $\mu_i:=\nu_i$ for all
$i \in  V(\Gamma_{P,i_o})$ and $\mu_i := 0$ otherwise. Therefore
$\mu_{i_o}=1$, which shows that $\underline{\mu}\neq 0$.
The equalities (\ref{eqlim}) and formula (\ref{formint}) imply:
$$Q_{(\Gamma, \underline{e}^{(k_o)})} (\underline {\mu}) =\,0$$
which contradicts the fact that $Q_{(\Gamma, \underline{e}^{(k_0)})}$ is
negative definite. \hfill $\Box$
\medskip

As an immediate consequence of Theorem \ref{mainth} we get:

\begin{corollary} \label{nonzero}\label{l1}
   If a vertex $i \in V(\Gamma)$ satisfies
   that $n_i\neq 0$ for every   solution of
   (\ref{eqgen}), then $e_i$ takes only a finite number of
   values. Therefore,
   the system (\ref{eqgen})  has a finite number of solutions
   $(\underline{n}, \underline{e})\in (\N^*)^{V(\Gamma)}\times
   (\N^*)^{V(\Gamma)}$. Thus, if all vertices satisfy  $n_i\neq 0$
   for every solution of
   (\ref{eqgen}), then there are a finite number of possible weights
   (self-intersections)
   for the vertices making the graph
   n-Gorenstein.
\end{corollary}

\noindent
{\bf Proof.} If $n_i\neq 0$ for all $i \in V(\Gamma)$, we get from
the equations (\ref{eqgen}) that:
  $$e_i= \frac{1}{n_i}( q_i + \sum_{\stackrel{j\in
      V(\Gamma)}{j\leftrightarrow i}} e_{ij}
    n_j ).$$
As by Theorem \ref{mainth} there are only finitely many possibilities
for $\underline{n}$, the conclusion follows.\hfill $\Box$
\medskip

\section{The possible non-finiteness of
   self-intersections}
\label{nonfin}

In this section we impose again the restrictions $q_i\geq -2$, satisfied
when the system
(\ref{eqgen}) corresponds to potential normal surface singularities
(see the relations
(\ref{change})).  We want to describe to what extent not only
the weights $\underline{n}$ can take a finite number of values
(which is ensured by Theorem \ref{mainth}),
but also  $\underline{e}$.

In view of Corollary \ref{l1}, let us concentrate our attention on the
vertices $i_0$ and on the solutions $(\underline{n}, \underline{e})$
of (\ref{eqgen}) such that $n_{i_0}=0$.

\begin{proposition} \label{l2}
   If there exists a vertex $i_0$ such that $n_{i_0}=0$ for some
   solution $(\underline{n}, \underline{e})$ of the system
  (\ref{eqgen}), then $v_{i_0}\leq 2$. Moreover:

 {\bf a)}  If  $v_{i_o} = 2$, then one
 has a cusp graph  and $n_i = 0$ for all $i\in V(\Gamma)$. Therefore $Z_K= \sum
 E_i$, $e_j\geq 2$ for all vertices $j$ and $e_j\geq 3$ for at least
 one vertex (in order to ensure that the graph is negative definite).

{\bf b)}  If  $v_{i_o} = 1$,
then $p_{i_0}=0$
   and the unique neighbor $i_1$ of $i_0$ satisfies $n_{i_1}=1$.
Therefore, $E_{i_o}$ is a smooth rational curve, $z_{i_0}=1$ and
$z_{i_1}=2$.

{\bf c)} If  $v_{i_0}=0$
(therefore $\Gamma$ has only one vertex and no edges), then
$p_{i_0}=1$. For every choice of weight $e_{i_0}$, we get a
n-Gorenstein graph and $Z_K= E_{i_0}$.
 \end{proposition}

 \noindent
{\bf Proof.} From (\ref{eqgen}) and (\ref{change}), we get:
  \begin{equation} \label{rel0}
   0=v_{i_0}+2p_{i_0}-2+\sum_{\stackrel{j \in
       V(\Gamma)}{j\leftrightarrow i}}
    e_{i_0j}n_j.
  \end{equation}
  We conclude using the fact that all the parameters appearing in the
  equation are non-negative integers. \hfill $\Box$
\medskip

Conversely, as shown by the examples iv) of Section \ref{exampl},
for all cusp-graph $\Gamma$ one has an infinite number of
possibilities for the weight $e_i$, and this for each vertex $i$.
More precisely, the set of possibilities for $\underline{e}$ is
$(\N^*\setminus \{1\})^{V(\Gamma)} \setminus (\{2\})^{V(\Gamma)}$.
The case b) of the previous lemma is realized for example in the
family iii) of Section \ref{exampl}. The case c) corresponds, for
instance, to the Brieskorn singularities defined by the equation
$z_1^a + z_2^b + z_3^c=0$ with $1/a + 1/b + 1/c = 1$;
these have a minimal resolution which
is a line bundle over a non-singular elliptic curve.

\medskip

The following  is an immediate consequence of Proposition \ref{l2}.

\begin{proposition} \label{measfin}
  Suppose that $(\Gamma, \underline{p})$ is not a cusp graph. We fix a
  vertex $i$ of   $\Gamma$ and
  we look at the possible values of $e_i$, when one varies $(\underline{n},
  \underline{e})$ among the solutions of (\ref{eqgen}). Then:

  $\bullet$ if $v_i>1$, the number $e_i$ takes only a finite
  number of values;

  $\bullet$ if  $v_i=1$, the number $e_i$ takes an infinite number of
  values if and
     only if there exists a solution with $n_i=0$. In this case, one
     obtains new
     solutions by varying only $e_i$ in $\N^*\setminus \{1\}$ and
     fixing all the
     other values of $(\underline{n}, \underline{e})$.
\end{proposition}

As a conclusion, we ask some questions that arise naturally
from our work:

$\bullet$ Given $(\Gamma, \underline{p})$, is the set of weights
$\underline{e}$ making $(\Gamma, \underline{p}, \underline{e})$
minimal and n-Gorenstein always non-empty ?

$\bullet$ Given $(\Gamma, \underline{p})$, can one give explicit
upper bounds on the values of the functions $\underline{n}$ which
can be extended to a solution of (\ref{eqgen}) ?

$\bullet$ Given $(\Gamma, \underline{p})$, can one give an algorithm
to compute the values of the functions $\underline{n}$ which
can be extended to a solution of (\ref{eqgen}) ?

\section*{Acknowledgments}
The research of the second author was supported by CONACYT and
DGAPA-UNAM, Mexico.

\bibliographystyle{plain}

\begin{thebibliography}{MMM}

\bibitem{Artin} M. Artin, On isolated rational singularities of
surfaces,  \textit{Amer. J. Math.}  \textbf{88}  (1966), 129-136.

\bibitem{BHPV 04} W. P. Barth, K. Hulek, C.A.M. Peters and A. Van de Ven,
   \textit{Compact complex surfaces},  Second enlarged edition,
   (Springer, 2004).

\bibitem{B 00} E. Brieskorn,  Singularities in the work of
    Friedrich Hirzebruch, in \textit{Surveys in Differential Geometry}
    vol. VII, 
  (International Press, 2000), pp.17-60.

\bibitem{Do 74}  I. V. Dolgachev, Quotient conical singularities on
complex surfaces, \textit{Funct Anal. Appl.} \textbf{8} (1974) 160-161.


\bibitem{D 83}  I. V. Dolgachev, On the link space of a Gorenstein
   quasi-homogeneous surface singularity, \textit{Math. Ann.} \textbf{265}
 (1983) 529-540.

\bibitem{D 78} A. Durfee, The signature of smoothings of
    complex surface   singularities, \textit{Math. Ann.} \textbf{232}
  (1978) 85-98. 

\bibitem{D 79} A. Durfee, Fifteen characterizations of
    rational double points and simple critical points, 
  \textit{Enseignem. Math.} \textbf{25} (1979) 131-163.

\bibitem{G 62} H. Grauert, {\"U}ber Modifikationen und exzeptionnelle
  analytische Mengen, \textit{Math. Ann.} \textbf{146} (1962) 331-368.

\bibitem{JLY 07} L. Jia, H.S. Luk and  St. S.-T. Yau, From CR
    geometry to 
  algebraic geometry and combinatorial geometry, \textit{AMS/IP
  Studies in Advanced 
  Mathematics}  \textbf{39} (2007) 109-124.

 \bibitem{LS 95}  F. Larri{\'o}n and  J. Seade, Complex surface
     singularities from
    the combinatorial point of view, \textit{Topology and its
  Appl.} \textbf{66}  (1995) 251-265.

 \bibitem{L 87} H. B. Laufer, The multiplicity of isolated
     two-dimensional
    hypersurface singularities, \textit{Trans. of the AMS} \textbf{302} no. 2
  (1987) 489-496.

\bibitem{L 84} E. Looijenga, \textit{Isolated singular points on
    complete intersections} (Cambridge Univ. Press, 1984).

\bibitem{M 61} D. Mumford, The topology of normal
    singularities of an algebraic surface and a criterion for
    simplicity, \textit{Publ. Math. IHES} \textbf{9} (1961) 229-246.


\bibitem{Ne 83} W. Neumann, The geometry of quasihomogeneous
surface singularities, in \textit{Proc. Symp. Pure Maths.} vol
\textbf{40} (A.M.S., 1983), pp.  245-258.


 \bibitem{R 97} M. Reid, Chapters on Algebraic Surfaces, in
  \textit{Complex Algebraic Geometry}, J. Koll{\'a}r editor, (A.M.S.,
  1997),  pp. 3-159.

  \bibitem{V 44} P. Du Val, On absolute and non-absolute
    singularities of algebraic surfaces, \textit{Revue de la Facult{\'e} des
    Sciences de l'Univ. d'Istanbul} (A) \textbf{91} (1944) 159-215.

 \bibitem{W 04}  C.T.C. Wall, \textit{Singular points of plane curves}
     (Cambridge Univ. Press, 2004).

 \bibitem{Y 79} St. S.-T. Yau, Hypersurface weighted dual
     graphs of normal 
   singularities of surfaces, \textit{Amer. Journ. of Maths.}
 \textbf{101} no.4 (1979)   761-812.



\end{thebibliography}

\end{document}